\documentclass[aip,graphicx,amsfonts,groupedaddress]{revtex4-1}

\draft 
\usepackage{tikz}
\usetikzlibrary{calc}
\newcommand{\ben}{\begin{equation}}
\newcommand{\bal}{\begin{align}}
\newcommand{\een}{\end{equation}}
\newcommand{\eal}{\end{align}}
\newcommand{\bea}{\begin{eqnarray}}
\newcommand{\eea}{\end{eqnarray}}

\def\Z{{\mathbb Z}} \def\R{{\mathbb R}}  \def\N{{\mathbb N}}

\begin{document}


\title{Polytope Expansion of Lie Characters and Applications} 



\author{Mark A. Walton\,}
\email[]{walton@uleth.ca}
\affiliation{Department of Physics and Astronomy, University of Lethbridge, Lethbridge, Alberta, Canada\ \ T1K 3M4}


\date{\today}

\begin{abstract}
 The weight systems of finite-dimensional representations of complex, simple Lie algebras exhibit patterns beyond Weyl-group symmetry. These patterns occur because weight systems can be decomposed into lattice polytopes in a natural way. Since lattice polytopes are relatively simple, this decomposition is useful, in addition to being more economical than the decomposition into single weights. An expansion of characters into polytope sums follows from the polytope decomposition of weight systems. We study this polytope expansion here. A new, general formula is given for the polytope sums involved. The combinatorics of the polytope expansion are analyzed;  we point out that they are reduced from  those of the Weyl character formula (described by the Kostant partition function) in an optimal way. We also show that the weight multiplicities can be found easily from the polytope multiplicities, indicating explicitly the equivalence of the two descriptions. Finally, we demonstrate the utility of the polytope expansion by showing how polytope multiplicities can be used in the calculation of tensor product decompositions, and subalgebra branching rules.
\end{abstract}


\maketitle 

\section{Introduction}

Lie groups and algebras have played an important role in physics, and the characters of their representations have been indispensable tools. Although much simpler than the full representations, characters can yield important physical information.  They encode physical spectra, and their algebra can describe couplings between different (sub-)systems.  In addition,  the decompositions of Lie characters into subgroup (or subalgebra)  characters can go a long way toward describing the physics of symmetry  breaking.

Here we study the characters of the integrable, unitary, finite-dimensional, highest-weight representations of the simple, complex Lie algebras. Specifically, we develop the general formula for these characters that we conjectured\cite{Walton} in 2002 and that was proved (independently) in Ref.~\onlinecite{Schutzer}. Some related results were recorded in Ref.~\onlinecite{SchutzerPhD}.

In Ref.~\onlinecite{Schutzer}, the formula was presented as a new character formula. In Ref.~\onlinecite{Walton}, on the other hand, it was treated as an expansion of the famous Weyl character formula motivated by the theory of lattice polytopes.\cite{B} Here we focus on  that interpretation of the formula, and develop it further.

The polytope expansion gives the character formula a geometric motivation. In addition, the formula displays patterns of weight multiplicities in an irreducible representation in as simple a manner as possible. The weight system of an irreducible representation is decomposed into multiplicity-1 subsets that form lattice polytopes. (For a simple illustration, see Figure \ref{figure:PolytopeExpn} below.)

The character is then a sum of the corresponding polytope generating functions (or sums). Here we derive a new, general formula for the polytope sums. Equation (\ref{Blan}) is a crucial component of the polytope expansion, and completes the description of the new character formula.\cite{Walton, Schutzer}

We show how weight multiplicities can be recovered easily from the polytope expansion (see the result (\ref{mA}) below). The equivalence of the weight and polytope systems is thus demonstrated explicitly.  We also argue that the combinatorics of the polytope expansion formula are an optimal improvement of those of the Weyl character formula.

The multiplicity-1 property of lattice polytopes makes the polytope sums relatively easy to use.  The polytope expansion formula therefore naturally leads to applications. We show here how it can be used to compute tensor product coefficients and subalgebra branching rules: equations (\ref{TAU},\,\ref{Udelta}) and (\ref{bApbAinv}) are the results. Just the polytope multiplicities appear in these formulas; the weight multiplicities are not necessary.  The polytope expansion not only displays patterns and is economical, it is useful in applications.

The following section discusses the Weyl character formula, to establish background and notation. Section 3 describes and develops the polytope expansion of Lie characters.  Section 4 first shows how weight multiplicities can be recovered simply from polytope multiplicities, and then treats two applications:  the calculation of  tensor product decompositions and subalgebra branching rules.

\section{Weyl character formula}

For later comparison and to establish notation, we will first discuss the Weyl character formula and related results that are relevant to the polytope expansion of Lie characters.

Let $P = \Z\{\Lambda^i\,|\, i=1,\ldots,r \}$
denote the weight lattice of a simple Lie algebra $X_r$, of rank
$r$. Here $\Lambda^i$ stands for the $i$-th fundamental weight.
$R_+$ ($R_-$) will represent the set of positive (negative) roots
of $X_r$, and
$S=\{\alpha_i\,|\, i=1,\ldots,r \}$ its simple roots.

The highest weights of integrable irreducible representations of
$X_r$ belong to the set $P_+ = \{\lambda = \sum_{j=1}^r\lambda_j\Lambda^j
\,|\, \lambda_j\in \N_{0} \}$, where $\N_0=\{0,1,2,\ldots\}$. For $\lambda,\mu\in P$, we write $\lambda\geq\mu$ if $\lambda-\mu\in \N_0 S$.

Let $W = \langle r_i\,|\, i=1,\ldots,r \rangle$ be the Weyl group of $X_r$,
where $r_i$ is the $i$-th
primitive Weyl reflection, with action $r_i\mu =
\mu - \mu_i\alpha_i$ on $\mu\in P$. $\det w$ is the sign of $w\in W$, and
$w.\lambda = w(\lambda+\rho)-\rho$ is the shifted action of $w$ on $\lambda$,
with the Weyl vector $\rho = \sum_{i=1}^r \Lambda^i =
\sum_{\alpha\in R_+} \alpha/2$.

The irreducible $X_r$-representation of highest weight $\lambda\in P_+$ will be denoted $L(\lambda)$. The (formal) character of $L(\lambda)$ can be expressed as
\begin{equation}
{\rm ch}_\lambda\ =\ \sum_{\mu\in P}\, {\rm mult}_\lambda(\mu)\, e^\mu\ ,
\label{chmult}\end{equation}
as an expansion in weights. Here ${\rm mult}_\lambda(\mu)$ denotes the multiplicity of weight $\mu$
in $L(\lambda)$, and $e^\mu$ is the formal exponential of the weight $\mu$. The weight system $P(\lambda)$ is the set of weights with non-zero multiplicities, \ben P(\lambda)\ =\ \{\, \mu\in P_+\ \vert\ {\rm mult}_\lambda(\mu)\in\N   \,\}  .\label{Pofla}\een
The sum in (\ref{chmult}) has been extended to $\mu\in P$ by setting ${\rm mult}_\lambda(\mu)=0$ for all $\mu\not\in P(\lambda)$.

A nice expression for the character is \begin{equation}
{\rm ch}_\lambda\ =\
\sum_{w\in W}\, e^{w\lambda}\, \prod_{\alpha\in R_+}\,
\lfloor 1-e^{-w\alpha}\rfloor^{-1}\ .\label{fchar}
\end{equation}
Here the Weyl invariance of the character is manifest. The price paid is that the expression is formal, emphasized here by the notation
\ben
\lfloor 1-e^\beta \rfloor^{-1}\ :=\ \left\{\matrix{(1-e^\beta)^{-1} = 1+e^\beta+e^{2\beta}+\ldots\ &,\ \beta\in R_-; \cr  -e^{-\beta}(1-e^{-\beta})^{-1}\ = -e^{-\beta}-e^{-2\beta}-\ldots &,\  \beta\in R_+. } \right. \label{formal}
\een
The basic rule-of-thumb is to expand in  terms of the type $e^\beta$, with $\beta$ a negative root.

\begin{figure}[ht]
  \begin{tikzpicture}[scale=0.5]
  \coordinate (Lone) at (0.8660254,0.5);
  \coordinate (Ltwo) at (0,1);
  \coordinate (aone) at (1.73205,0);
  \coordinate (atwo) at (-0.8660254,1.5);
  \coordinate (aonetwo) at ($(aone)+(atwo)$);
  \coordinate (eps) at (0.15,0.15);
  \coordinate (la) at ($2*(Lone)+5*(Ltwo)$);
  \coordinate (abb0) at ($(la)-6.5*(aone)$);
  \coordinate (ab0) at ($(la)-6.*(aone)$);
  \coordinate (abs) at ($(la)-6.*(aone)-0.15*(atwo)$);
  \coordinate (ab2) at ($(la)-6.*(aone)-2*(atwo)$);
  \coordinate (a2s) at ($(la)-2*(aone)-0.15*(atwo)$);
  \coordinate (a2b) at ($(la)-2*(aone)-7.5*(atwo)$);
  \coordinate (a0bb) at ($(la)-7.5*(atwo)$);
  \coordinate (ab4) at ($(la)-8*(aone)-4*(atwo)$);
  \coordinate (as4) at ($(la)-0.15*(aone)-4*(atwo)$);
  \coordinate (asb) at ($(la)-0.15*(aone)-8*(atwo)$);
  \coordinate (a2s6) at ($(la)-1.85*(aone)-6*(atwo)$);
  \coordinate (a2sb) at ($(la)-1.85*(aone)-8*(atwo)$);
  \coordinate (abs6) at ($(la)-9*(aone)-6*(atwo)$);
  \tikzstyle{every node}=[draw,shape=circle, scale=0.75];
  \path (la) node (a00) {1};
  \path ($(la)-(aone)$) node (a10) {1};
  \path ($(la)-(atwo)$) node (a01) {1};
  \path ($(la)-(atwo)-(aone)$) node (a11) {2};
  \path ($(la)-(atwo)-2*(aone)$) node (a21) {1};
  \path ($(la)-2*(atwo)$) node (a02) {1};
  \path ($(la)-2*(atwo)-(aone)$) node (a12) {2};
  \path ($(la)-2*(atwo)-2*(aone)$) node (a22) {2};
  \path ($(la)-2*(atwo)-3*(aone)$) node (a32) {1};
  \path ($(la)-3*(atwo)$) node (a03) {1};
  \path ($(la)-3*(atwo)-(aone)$) node (a13) {2};
  \path ($(la)-3*(atwo)-2*(aone)$) node (a23) {2};
  \path ($(la)-3*(atwo)-3*(aone)$) node (a33) {2};
  \path ($(la)-3*(atwo)-4*(aone)$) node (a43) {1};
  \path ($(la)-4*(atwo)-(aone)$) node (a14) {1};
  \path ($(la)-4*(atwo)-2*(aone)$) node (a24) {1};
  \path ($(la)-4*(atwo)-3*(aone)$) node (a34) {1};
  \path ($(la)-4*(atwo)-4*(aone)$) node (a44) {1};
  \draw [<-, ultra thick] (abb0) -- (a10);
  \draw [ultra thick] (a10) -- (a00) -- (a01) -- (a02) -- (a02) -- (a03);
  \draw [->, ultra thick] (a03) -- (a0bb);
  \draw [<->, ultra thick] ($-10.*(Lone)+5*(Ltwo)+0.15*(atwo)-0.5*(aone)$) -- ($-6*(Lone)+3*(Ltwo)+0.15*(atwo)$) -- ($-2*(Lone)-5*(Ltwo)$);
  \draw [<->, ultra thick] (abs6) -- (a2s6) -- (a2sb);
  \draw [<-, ultra thick, dashed] (abs) -- (a2s) -- (a21);
  \draw [-, ultra thick, dashed] (a21) -- (a22) -- (a23) -- (a24);
  \draw [->, ultra thick, dashed] (a24) -- (a2b);
  \draw [<-, ultra thick, dashed] (ab4) -- (a44);
  \draw [-, ultra thick, dashed] (a44) -- (a34) -- (a24) -- (a14);
  \draw [->, ultra thick, dashed] (a14) -- (as4) -- (asb);
  \draw [<->, ultra thick, dashed] ($-4*(Lone)-1*(Ltwo)-1.5*(eps)-3.2*(aone)$) -- ($-4*(Lone)-1*(Ltwo)-1.5*(eps)+0.05*(aone)$) -- ($-4*(Lone)-1*(Ltwo)-1.5*(eps)-2.2*(atwo)+0.05*(aone)$);
  \end{tikzpicture}
  \caption{Weight diagram for $A_2$ representation of highest weight $\lambda=\Lambda^1+3\Lambda^2$, illustrating the Weyl character formula. Compare with Figure \ref{figure:Brion}. The numbers indicate weight multiplicities. Cones are drawn as 2 arrows emanating from a common vertex in the directions of the negative simple roots.  Each  represents a term in the sum of (\ref{WCF}) over elements of the Weyl group $W$. Negative (positive) terms, for $\det w =-1\, (+1)$,  are indicated by dashed (solid) lines. The multiplicities in the expansions of each such term (not indicated) are determined by the Kostant partition function (\ref{Kostant}). }
  \label{figure:WCF}
\end{figure}
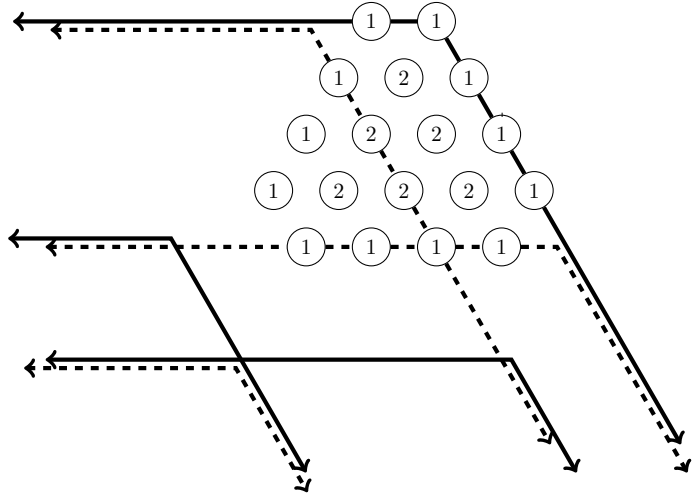

The usual formula can be recovered as follows. Each
Weyl element $w\in W$ separates the
positive roots into two disjoint sets:
\begin{equation}
\matrix{
R_+^w:= \{\alpha\in R_+\,|\, w\alpha\in R_+\}\ ,&
\ \ R_-^w:=
\{\alpha\in R_+\,|\, w\alpha\in R_-\}\ ,\cr R_+^w \cup R_-^w\ =\
R_+\ ,&
\ \ R_+^w \cap R_-^w\ =\ \{\}\ ,\cr wR_+^w\ =\ R_+^w\ \ \ \
,&
\ \ \ \ wR_-^w\ =\ -R^w_-\ \  .
}
\label{Rw}
\end{equation}
It can be shown that $\det w = (-1)^{\dim R^w_-}$, and
\begin{equation}
w\rho -\rho\ =\ -\sum_{\gamma\in R^w_-}\, \gamma\ =\ \sum_{\gamma\in R^w_-}\, w\gamma\ \ .
\end{equation}
Using these results, (\ref{fchar}) becomes
\bea
{\rm ch}_\lambda\ =\ \sum_{w\in W}\, e^{w\lambda}\,
\prod_{\beta\in R^w_-}\, (-e^{w\beta})(1-e^{w\beta})^{-1}\,
\prod_{\alpha\in R^w_+}\, (1-e^{-w\alpha})^{-1}\ \cr
\quad =\, \sum_{w\in W}\,(\det w)\, e^{w\lambda-w\sum_{\gamma\in R^w_-}\gamma}\,
\prod_{\beta\in R^w_-}\, (1-e^{w\beta})^{-1}\,
\prod_{\alpha\in R^w_+}\, (1-e^{-w\alpha})^{-1}\ .
\eea
The famous Weyl character formula
\begin{equation}
{\rm ch}_\lambda\ =\
\frac{\sum_{w\in W}\, (\det w)\, e^{w.\lambda}}
{\prod_{\alpha\in R_+}\, (1-e^{-\alpha})}\ \label{WCF}
\end{equation} then follows. See Figure \ref{figure:WCF} for an illustration using a representation of the rank-2 algebra $A_2$.

After using ${\rm ch}_0 = 1$ to derive the Weyl denominator formula, we can
rewrite the character formula as
\begin{equation}
{\rm ch}_\lambda\ =\
\frac{\sum_{w\in W}\, (\det w)\, e^{w(\lambda+\rho)}}
{\sum_{w\in W}\, (\det w)\, e^{w\rho}}\ . \label{WCFw}
\end{equation}
Comparing to (\ref{chmult}) confirms the Weyl symmetry
\begin{equation}
{\rm mult}_\lambda(\mu)\ =\ {\rm mult}_\lambda(w\mu)\ \ \ \ (\forall w\in W)
\label{Wsymm}\end{equation}
of the weight multiplicities.

Using the Weyl character formula to extend the definition of ${\rm ch}_\lambda$ to all $\lambda\in P$, we find another symmetry: \ben {\rm ch}_{w.\lambda}\ =\ (\det w)\, {\rm ch}_\lambda\ \ ,\ \ \ \forall w\in W\ .  \label{chW}\een

Let $\theta$ denote the highest root of $X_r$, the highest weight of its adjoint representation. Many multiplicities ${\rm mult}_\lambda(\mu)$  become very large as the  height $\langle\lambda,\theta\rangle$ of the highest weight increases  (except for rank $r=1$). The denominator of (\ref{WCF}) is independent of the highest weight $\lambda$ and so provides a measure of this rapid increase that is independent of the representation. The Kostant partition function $K(\gamma)$, $\gamma\in \N_0R_+$, is defined from that denominator by
\ben  \prod_{\alpha\in R_+}\, (1-e^{\alpha})^{-1}\ =\ \sum_{\gamma\in \N_0R_+}\, K(\gamma)\, e^\gamma\ .\label{Kostant}\een $K(\gamma)$ is the number of ways $\gamma$ can be written as a non-negative integer linear combination of positive roots.

The Kostant multiplicity formula \ben {\rm mult}_\lambda(\mu)\ =\ \sum_{w\in W}\, (\det w)\, K(w.\lambda-\mu) \label{KMF}\een can be derived by comparing (\ref{chmult}) and  (\ref{WCF}), using (\ref{chW}). Multiplicities increase rapidly with the height of the  highest weight  because of the number of positive roots. One important consideration they must take into account is the Weyl degeneracy (\ref{Wsymm}).

On the other hand, it is well known that the pattern of weight multiplicities makes plain degeneracies beyond the Weyl symmetry. We will now write a character formula \cite{Walton, Schutzer} that makes this further symmetry clear. We will be directed by the theory of lattice polytopes.\cite{B}

\section{Polytope sums and character expansion}

A polytope is the convex hull of finitely many points in $\R^d$
(see Ref.~\onlinecite{Zi}, for example). A lattice polytope (also known as an integral
polytope) is a polytope with all
its vertices in an integral lattice $\Lambda\in \R^d$
(see Ref.~\onlinecite{B}, for example). The
exponential sum of a lattice polytope $Pt$,
\begin{equation}
\sum_{x\in Pt\cap \Lambda}\, \exp\{\langle c, x\rangle\}\ ,
\end{equation}
is a useful tool. Here $c$ is a vector in $\R^d$, and
$\langle\cdot,\cdot\rangle$ is the usual inner product.
For comparison with formal characters (and simplicity of notation), we'll
consider formal exponential sums
\begin{equation}
E[Pt; \Lambda]\ :=\ \sum_{x\in Pt\cap\Lambda }\, e^x\ .
\end{equation}
These exponential sums are the
generating functions for the integral points in a lattice polytope.

Brion \cite{Br} has proved a formula for the exponential sum of a convex
lattice polytope, that expresses it as a sum of simpler terms,
associated with each of the vertices of the polytope. It reads
\begin{equation}
E[Pt; \Lambda ]\ =\ \sum_{x\in Pt\cap\Lambda}\, e^x\ =\
\sum_{v\in {\rm Vert}Pt}\, e^v\, V_v\ .
\end{equation}
Here ${\rm Vert}Pt$ is the set of vertices of $Pt$. The vertex term $V_v$ is the exponential sum of the cone of feasible directions at the vertex $v$.\cite{Br}

The form (\ref{WCF}) of the Weyl character formula is
similar to the Brion formula. To make this more precise, we'll write
the Brion formula for the exponential sum of
the lattice polytope $Pt_\lambda$ with vertices in the
Weyl orbit of a highest weight $\lambda\in P_\ge$, i.e.
${\rm Vert}Pt_\lambda =
W\lambda$. Let us call $Pt_\lambda$ the
polytope of highest weight $\lambda$.
The appropriate lattice here is the root lattice $Q$ of $X_r$,
shifted by $\lambda$: $\lambda + Q\subset P$:
\begin{equation}
Pt_\lambda\cap (\lambda+Q)\ =\ P(\lambda)\ .
\end{equation}
The polytope sum will be
\begin{equation}
E[Pt_\lambda; \lambda+Q]\ =\ \sum_{\mu\in P(\lambda)}\, e^\mu\ .
\end{equation}

As noted in Ref.~\onlinecite{Walton}, the Brion formula gives
\begin{equation}
B_\lambda\ =\
\sum_{w\in W}\, e^{w\lambda}\, \prod_{\alpha\in S}\,
\lfloor 1-e^{-w\alpha}\rfloor^{-1}\ ,
\label{Blaf}\end{equation}
for the weight polytope $Pt_\lambda$.

\begin{figure}[ht]
  \begin{tikzpicture}[scale=0.45]
  \coordinate (shift) at (12,0);
  \coordinate (Lone) at (0.8660254,0.5);
  \coordinate (Ltwo) at (0,1);
  \coordinate (aone) at (1.73205,0);
  \coordinate (atwo) at (-0.8660254,1.5);
  \coordinate (aonetwo) at ($(aone)+(atwo)$);
  \coordinate (bit) at (2,0);
  \path ($-2*(shift)+5*(Lone)-2*(Ltwo)+(bit)$) node {$\hookrightarrow$};
  \path ($-0.95*(shift)+5*(Lone)-2*(Ltwo)+0.625*(bit)$) node {$+$};
  \path ($-0.175*(shift)+3*(Lone)-1*(Ltwo)+0.625*(bit)$) node {$+$};
  \tikzstyle{every node}=[draw,shape=circle, scale=0.7];
  %
  %
  \path ($-2*(shift)+3*(Lone)+2*(Ltwo)$) node {1};
  \path ($-2*(shift)+1*(Lone)+3*(Ltwo)$) node {1};
  \path ($-2*(shift)-1*(Lone)+4*(Ltwo)$) node {1};
  \path ($-2*(shift)-3*(Lone)+5*(Ltwo)$) node {1};
  \path ($-2*(shift)-4*(Lone)+4*(Ltwo)$) node {1};
  \path ($-2*(shift)-5*(Lone)+3*(Ltwo)$) node {1};
  \path ($-2*(shift)-4*(Lone)+1*(Ltwo)$) node {1};
  \path ($-2*(shift)-3*(Lone)-1*(Ltwo)$) node {1};
  \path ($-2*(shift)-2*(Lone)-3*(Ltwo)$) node {1};
  \path ($-2*(shift)+0*(Lone)-4*(Ltwo)$) node {1};
  \path ($-2*(shift)+2*(Lone)-5*(Ltwo)$) node {1};
  \path ($-2*(shift)+3*(Lone)-4*(Ltwo)$) node {1};
  \path ($-2*(shift)+4*(Lone)-3*(Ltwo)$) node {1};
  \path ($-2*(shift)+5*(Lone)-2*(Ltwo)$) node {1};
  \path ($-2*(shift)+4*(Lone)+0*(Ltwo)$) node {1};
  \path ($-2*(shift)+2*(Lone)+1*(Ltwo)$) node {2};
  \path ($-2*(shift)+0*(Lone)+2*(Ltwo)$) node {2};
  \path ($-2*(shift)-2*(Lone)+3*(Ltwo)$) node {2};
  \path ($-2*(shift)-3*(Lone)+2*(Ltwo)$) node {2};
  \path ($-2*(shift)-2*(Lone)+0*(Ltwo)$) node {2};
  \path ($-2*(shift)-1*(Lone)-2*(Ltwo)$) node {2};
  \path ($-2*(shift)+1*(Lone)-3*(Ltwo)$) node {2};
  \path ($-2*(shift)+2*(Lone)-2*(Ltwo)$) node {2};
  \path ($-2*(shift)+3*(Lone)-1*(Ltwo)$) node {2};
  \path ($-2*(shift)+1*(Lone)+0*(Ltwo)$) node {3};
  \path ($-2*(shift)-1*(Lone)+1*(Ltwo)$) node {3};
  \path ($-2*(shift)+0*(Lone)-1*(Ltwo)$) node {3};
  %
  %
  \path ($-0.95*(shift)+3*(Lone)+2*(Ltwo)$) node {1};
  \path ($-0.95*(shift)+1*(Lone)+3*(Ltwo)$) node {1};
  \path ($-0.95*(shift)-1*(Lone)+4*(Ltwo)$) node {1};
  \path ($-0.95*(shift)-3*(Lone)+5*(Ltwo)$) node {1};
  \path ($-0.95*(shift)-4*(Lone)+4*(Ltwo)$) node {1};
  \path ($-0.95*(shift)-5*(Lone)+3*(Ltwo)$) node {1};
  \path ($-0.95*(shift)-4*(Lone)+1*(Ltwo)$) node {1};
  \path ($-0.95*(shift)-3*(Lone)-1*(Ltwo)$) node {1};
  \path ($-0.95*(shift)-2*(Lone)-3*(Ltwo)$) node {1};
  \path ($-0.95*(shift)+0*(Lone)-4*(Ltwo)$) node {1};
  \path ($-0.95*(shift)+2*(Lone)-5*(Ltwo)$) node {1};
  \path ($-0.95*(shift)+3*(Lone)-4*(Ltwo)$) node {1};
  \path ($-0.95*(shift)+4*(Lone)-3*(Ltwo)$) node {1};
  \path ($-0.95*(shift)+5*(Lone)-2*(Ltwo)$) node {1};
  \path ($-0.95*(shift)+4*(Lone)+0*(Ltwo)$) node {1};
  \path ($-0.95*(shift)+2*(Lone)+1*(Ltwo)$) node {1};
  \path ($-0.95*(shift)+0*(Lone)+2*(Ltwo)$) node {1};
  \path ($-0.95*(shift)-2*(Lone)+3*(Ltwo)$) node {1};
  \path ($-0.95*(shift)-3*(Lone)+2*(Ltwo)$) node {1};
  \path ($-0.95*(shift)-2*(Lone)+0*(Ltwo)$) node {1};
  \path ($-0.95*(shift)-1*(Lone)-2*(Ltwo)$) node {1};
  \path ($-0.95*(shift)+1*(Lone)-3*(Ltwo)$) node {1};
  \path ($-0.95*(shift)+2*(Lone)-2*(Ltwo)$) node {1};
  \path ($-0.95*(shift)+3*(Lone)-1*(Ltwo)$) node {1};
  \path ($-0.95*(shift)+1*(Lone)+0*(Ltwo)$) node {1};
  \path ($-0.95*(shift)-1*(Lone)+1*(Ltwo)$) node {1};
  \path ($-0.95*(shift)+0*(Lone)-1*(Ltwo)$) node {1};
  %
  %
  \path ($-0.175*(shift)+2*(Lone)+1*(Ltwo)$) node {1};
  \path ($-0.175*(shift)+0*(Lone)+2*(Ltwo)$) node {1};
  \path ($-0.175*(shift)-2*(Lone)+3*(Ltwo)$) node {1};
  \path ($-0.175*(shift)-3*(Lone)+2*(Ltwo)$) node {1};
  \path ($-0.175*(shift)-2*(Lone)+0*(Ltwo)$) node {1};
  \path ($-0.175*(shift)-1*(Lone)-2*(Ltwo)$) node {1};
  \path ($-0.175*(shift)+1*(Lone)-3*(Ltwo)$) node {1};
  \path ($-0.175*(shift)+2*(Lone)-2*(Ltwo)$) node {1};
  \path ($-0.175*(shift)+3*(Lone)-1*(Ltwo)$) node {1};
  \path ($-0.175*(shift)+1*(Lone)+0*(Ltwo)$) node {1};
  \path ($-0.175*(shift)-1*(Lone)+1*(Ltwo)$) node {1};
  \path ($-0.175*(shift)+0*(Lone)-1*(Ltwo)$) node {1};
  %
  %
  \path ($0.3125*(shift)+1*(Lone)+0*(Ltwo)$) node {1};
  \path ($0.3125*(shift)-1*(Lone)+1*(Ltwo)$) node {1};
  \path ($0.3125*(shift)+0*(Lone)-1*(Ltwo)$) node {1};
  \end{tikzpicture}
  \caption{Weight diagram for $A_2$ representation of highest weight $\lambda=3\Lambda^1+2\Lambda^2$, illustrating its polytope decomposition. For this example, the polytope expansion (\ref{chAB}) reads ${\rm ch}_{3\Lambda^1+2\Lambda^2} = B_{3\Lambda^1+2\Lambda^2} + B_{2\Lambda^1+\Lambda^2} + B_{\Lambda^1}$. Incidentally, for an irreducible $A_2$ representation of arbitrary highest weight $\lambda=\lambda_1\Lambda^1+\lambda_2\Lambda^2$, the expansion is ${\rm ch}_\lambda = B_\lambda+ B_{\lambda-\theta}+\ldots +B_{\lambda-{\rm min}(\lambda_1,\lambda_2)\,\theta}$\,, an elegant presentation of the well-known pattern of $A_2$ weight multiplicities.}
  \label{figure:PolytopeExpn}
\end{figure}
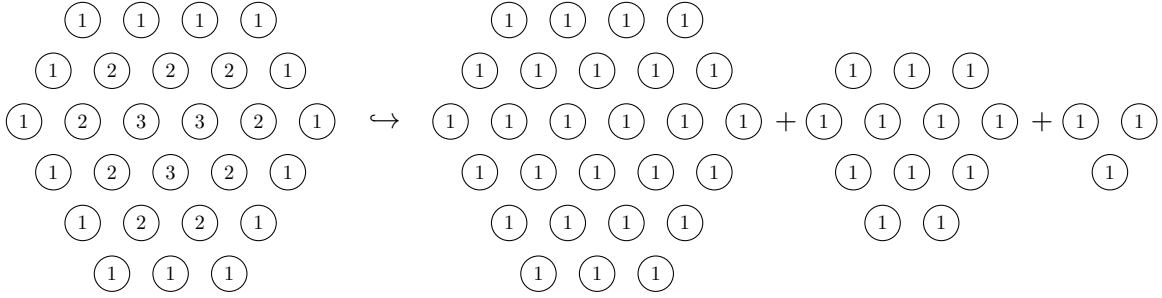

The similarity between the formal version (\ref{fchar}) of the Weyl character formula and the Brion formula (\ref{Blaf})
is striking.  Because lattice-polytope sums are simple, with cone sums (in $\{0,1\}\,e^P$) instead of partition functions (in $\N_0\, e^P$), it is sensible to exploit the similarity and write \ben  {\rm ch}_\lambda\ =\ \sum_{\mu\leq\lambda}\, A_{\lambda,\mu}\, B_\mu\ .\label{chAB}\een For an $A_2$ illustration of this polytope expansion, see Figure \ref{figure:PolytopeExpn}.

This polytope-expansion formula for Lie characters was reported in Ref.~\onlinecite{Walton}. We'll call the coefficients {\it polytope multiplicities}, and sometimes write  \ben A_{\lambda,\mu}\ =:\ {\rm polyt}_\lambda(\mu)\ , \label{polytA}\een in analogy with weight multiplicities ${\rm mult}_\lambda(\mu)$. They were conjectured to be non-negative integers in Ref.~\onlinecite{Walton}. Proofs of this and (\ref{chAB}) have been given recently by Sch\"utzer.\cite{Schutzer}

The $A_{\lambda,\mu}$ can be considered the entries of a matrix that is triangular, and so easily invertible. Furthermore, writing \ben B_\lambda\ =\ \sum_{\mu\leq\lambda}\, A^{-1}_{\lambda,\mu}\, {\rm ch}_\mu\  \label{BAinvch}\een leads to an explicit formula \cite{Walton} \ben A^{-1}_{\lambda,\mu}\ =\ \sum_{w\in W}\, (\det w)\, F(\lambda-w.\mu)\ .  \label{AinvF}\een  Here $F$ is defined by \begin{equation}
\prod_{\gamma\in R_+\backslash S}\, (1-e^{\gamma})\ =:\
\sum_{\beta\in \N_0 R_+}\, F(\beta)\, e^\beta\ ,\label{Fdef}
\end{equation} so that \begin{equation}
B_\lambda\ =\ \sum_{\beta\in \N_{0} R_+}\, F(\beta)\,
{\rm ch}_{\lambda-\beta}\ ,
\label{BFch}\end{equation} and (\ref{AinvF}) follows from (\ref{chW}).

So the polytope multiplicities can be found rather easily. First, the matrix with entries (\ref{AinvF}) is calculated in simple fashion, then inverted.  Again, the inversion is straightforward, since the matrix is triangular.   In Ref.~\onlinecite{Walton}, many polytope multiplicities were calculated for low-rank examples. Formulas for ``polytope dimensions'' (or discrete volumes of polytopes) were useful there, and follow from (\ref{AinvF}): \ben  b_\lambda\ =\ \sum_{\mu\leq \lambda}\, A^{-1}_{\lambda,\mu}\ d_\mu\ ,\label{blAinvdl}\een when the Weyl dimension formula \ben d_\lambda\ :=\ \dim\, L(\lambda)\ =\ \prod_{\alpha\in R_+}\, \langle \lambda+\rho,\alpha\rangle\, /\, \langle \rho,\alpha\rangle\  \label{Wdf}\een is used.

The relation (\ref{BFch}) can also be interpreted as a recursion relation for the polytope multiplicities $A_{\lambda,\mu}={\rm polyt}_\lambda(\mu)$.    Examples are given in Refs.~\onlinecite{Walton, Schutzer}, and the general relation \ben {\rm ch}_\mu\ =\ B_\mu\, -\, \sum_{\gamma<\mu}\, A^{-1}_{\mu,\gamma}\,{\rm ch}_\gamma\ =\ B_\mu\, -\, \sum_{\beta\not=0}\, F(\beta)\,
{\rm ch}_{\lambda-\beta}\ \label{recurs}\een is shown in Ref.~\onlinecite{Walton}. At least for simple low-rank cases, polytope multiplicities can be found recursively using (\ref{recurs}).

\begin{figure}[ht]
  \begin{tikzpicture}[scale=0.5]
  \coordinate (Lone) at (0.8660254,0.5);
  \coordinate (Ltwo) at (0,1);
  \coordinate (aone) at (1.73205,0);
  \coordinate (atwo) at (-0.8660254,1.5);
  \coordinate (aonetwo) at ($(aone)+(atwo)$);
  \coordinate (eps) at (0.15,0.15);
  \coordinate (la) at ($2*(Lone)+5*(Ltwo)$);
  \tikzstyle{every node}=[draw,shape=circle, scale=0.75];
  \path (la) node (a00) {1};
  \path ($(la)-(aone)$) node (a10) {1};
  \coordinate (a60) at ($(la)-6.5*(aone)$);
  \path ($(la)-(atwo)$) node (a01) {1};
  \path ($(la)-(atwo)-(aone)$) node (a11) {1};
  \path ($(la)-(atwo)-2*(aone)$) node (a21) {1};
  \coordinate (a07) at ($(la)-7.5*(atwo)$);
  \path ($(la)-2*(atwo)$) node (a02) {1};
  \path ($(la)-2*(atwo)-(aone)$) node (a12) {1};
  \path ($(la)-2*(atwo)-2*(aone)$) node (a22) {1};
  \path ($(la)-2*(atwo)-3*(aone)$) node (a32) {1};
  \path ($(la)-3*(atwo)$) node (a03) {1};
  \path ($(la)-3*(atwo)-(aone)$) node (a13) {1};
  \path ($(la)-3*(atwo)-2*(aone)$) node (a23) {1};
  \path ($(la)-3*(atwo)-3*(aone)$) node (a33) {1};
  \path ($(la)-3*(atwo)-4*(aone)$) node (a43) {1};
  \path ($(la)-4*(atwo)-(aone)$) node (a14) {1};
  \path ($(la)-4*(atwo)-2*(aone)$) node (a24) {1};
  \path ($(la)-4*(atwo)-3*(aone)$) node (a34) {1};
  \path ($(la)-4*(atwo)-4*(aone)$) node (a44) {1};
  \draw [<-, ultra thick] (a60) -- (a10);
  \draw [ultra thick] (a10) -- (a00) -- (a01) -- (a02) -- (a02) -- (a03);
  \draw [->, ultra thick] (a03) -- (a07);
  \draw [<->, ultra thick, dashed] ($-2*(Lone)+7*(Ltwo)-(eps)-4.02*(aone)-0.02*(atwo)$) -- ($-2*(Lone)+7*(Ltwo)-(eps)-0.02*(aone)-0.02*(atwo)$) -- ($-2*(Lone)+7*(Ltwo)-(eps)-7*(aonetwo)$);
  \draw [<->, ultra thick,dashed] ($6*(Lone)-3*(Ltwo)-(eps)-6.5*(aonetwo)$) -- ($6*(Lone)-3*(Ltwo)-(eps)$) -- ($6*(Lone)-3*(Ltwo)-(eps)-3*(atwo)$);
  \draw [<->, ultra thick,dashed] ($-4*(Lone)+2*(Ltwo)-4*(aonetwo)$) -- ($-4*(Lone)+2*(Ltwo)$) -- ($-4*(Lone)+2*(Ltwo)-6.5*(atwo)$);
  \draw [<->, ultra thick, dashed] ($3*(Lone)-3*(Ltwo)-7*(aone)$) -- ($3*(Lone)-3*(Ltwo)$) -- ($3*(Lone)-3*(Ltwo)-5*(aonetwo)$);
  \draw [<->, ultra thick] ($-3*(Lone)-(eps)-3.75*(aone)-0.03*(atwo)$) -- ($-3*(Lone)-(eps)-0.03*(atwo)$) -- ($-3*(Lone)-(eps)-5*(atwo)$);
  \end{tikzpicture}
  \caption{Lattice polytope in $A_2$ weight space $Pt_\lambda$ of highest weight $\lambda=\Lambda^1+3\Lambda^2$, illustrating the formal Brion formula  (\ref{Blaf}) and the more explicit version (\ref{Blan}). Compare with Figure \ref{figure:WCF}.  All ``weight multiplicities'' are 1, by definition. Each cone, drawn as 2 arrows emanating from a common vertex, represents a term in the sum of (\ref{Blan}) over elements of the Weyl group $W$. Negative (positive) terms are indicated by dashed (solid) lines. Up to sign, all ``cone multiplicities'' are 1.}
  \label{figure:Brion}
\end{figure}
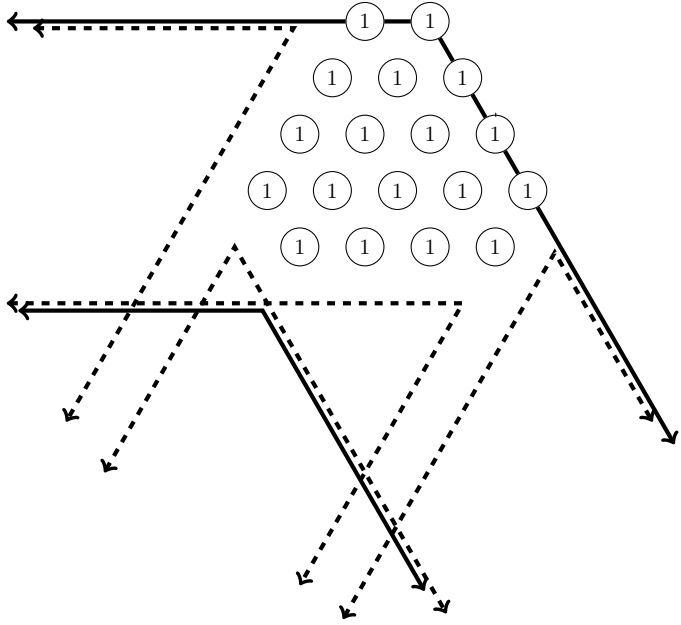

Let us now use (\ref{formal}) and rewrite equation (\ref{Blaf}). Define \ben wS\cap R_+\ =:\ S_+^w\ ,\ \ wS\cap R_-\ =:\ S_-^w\  . \label{Swpm}\een We get
\ben  B_\lambda\ =\ \sum_{w\in W}\, e^{w\lambda}\, \prod_{\beta\in S_+^w}\, (1-e^{-\beta})^{-1}\,  \prod_{\gamma\in S_-^w}\, (-e^\gamma)(1-e^\gamma)^{-1}\ . \label{Blai}\een
Define \ben \sigma(w)\ :=\ -\,\sum_{\gamma\in S^w_-}\, \gamma\ ,\ \ {\rm and}\ \ \epsilon(w)\ :=\ (-1)^{\dim S^w_-}\ \ .\label{siwmin}\een Notice that $\sigma(w)\in \N_0 R_+$ except that $\sigma(1)=0$. With \ben \vert w\alpha\vert\ :=\ \left\{\matrix{ w\alpha\, ,\ \ & w\alpha\in R_+\ ;\cr
-w\alpha\, ,\ \, & w\alpha\in R_-\ ,  }\right.  \label{absv}\een we can write
\ben  B_\lambda\ =\ \sum_{w\in W}\, e^{w\lambda}\ \epsilon(w)\,e^{- \sigma(w)}\, \prod_{\alpha\in S}\, (1-e^{-\vert w\alpha \vert})^{-1}\ .    \label{Blan}\een  Figure \ref{figure:Brion} illustrates the polytope formula of equations (\ref{Blan}) and (\ref{Blaf}).

Comparing with the Weyl character formula of equation (\ref{WCF}), we note an important difference. While the Kostant factor $\prod_{\alpha\in R_+}\, (1-e^{-\alpha})^{-1}$ is independent of $w\in W$, the polytope expansion  involves the factor $\prod_{\alpha\in S}\, (1-e^{-\vert w\alpha\vert})^{-1}$, a different term for each $w$ in the sum over $W$.  As indicated above, the compensation for this complication is trivial combinatorics.

To see this, generalize the Kostant partition function of (\ref{Kostant}) to \ben  \prod_{\alpha\in U}\, (1-e^{\alpha})^{-1}\ =\ \sum_{\gamma\in \N_0R_+}\, K_U(\gamma)\, e^\gamma\ ,  \label{gKostant}\een for any $U\subset R_+$. The original Kostant partition function is  $K=K_{R_+}$. Then if \ben \vert wS\vert\ :=\ \{\, \vert w\alpha\vert\ :\  \alpha\in S\, \}\ ,\label{vwSv}\een then $K_{\vert wS\vert}$ is the relevant partition function for (\ref{Blan}).  $K(\gamma)$ in (\ref{Kostant}) can become very large. In stark contrast,  since $\dim\vert wS\vert= \dim S = r$, the number of simple roots, \ben  K_{\vert wS\vert}(\gamma)\ =\ \left\{ \matrix{ 1\ ,\ &\ {\rm if\ } \gamma\in\N_0\, \vert wS\vert\ ;\cr 0\ ,\ &\ {\rm otherwise\, .}  }   \right.   \label{KwS}\een
This is consistent with the identification \ben V_{w\lambda}\ =\ \prod_{\alpha\in S}\, \lfloor 1-e^{-w\alpha}\rfloor^{-1}\ =\ \epsilon(w)\, e^{-\sigma(w)}\, \prod_{\alpha\in S}\, (1-e^{-\vert w\alpha \vert})^{-1}\ . \label{vertKw}\een
Notice that the expression just written is independent of $\lambda$. We can therefore write \ben V_{w\lambda}\ =:\ V_{w}\ ,\ \ \lambda\in P_+\ ,\label{Vw}\een and \ben B_\lambda\ =\ \sum_{w\in W}\, e^{w\lambda}\, V_{w}\ . \label{Vwla}\een

To close this section, we argue that the polytope expansion formula\cite{Walton,Schutzer} (\ref{chAB}) improves on the Kostant combinatorics of the Weyl character formula (\ref{WCF}) in optimal fashion. In the polytope expansion formula (\ref{chAB}), the Kostant partition function (\ref{Kostant}) of the Weyl charcater formula (\ref{WCF}) is replaced by a product of 2 simpler partition  functions. In terms of the subsets of positive roots that define the partition functions,  $R_+$ is replaced by $S\cup\left( R_+\backslash S \right)$. This is optimal in the sense that splitting off a bigger subset would introduce 2 non-trivial partition functions, instead of a single partition function and multiplicity-1 vertex (cone) generating functions.

Put another way, the lattice-polyhedral-cone expansion of the Kostant partition function is at the base of the polytope-expansion formula for Lie characters.  The relation
\ben  \sum_{\gamma\in \N_0R_+}\, K(\gamma)\, e^\gamma\ =\ V_{w=1}\ \prod_{\beta\in R_+\backslash S}\, (1-e^{\beta})^{-1}\, \  .  \label{Kcone}\een follows from (\ref{Kostant}) and (\ref{Vw}), for example.

\section{Weights from polytopes and applications of the polytope expansion}\label{sfour}

While characters involve the highly non-trivial weight multiplicities, polytope sums  are multiplicity-1.  Thus the polytope expansion (\ref{chAB}) of a Lie character ${\rm ch}_\lambda$ allows one to work with simpler objects $B_\mu$, $\mu\leq\lambda$. That is the main advantage in applications of the new character formula, as we show in this section.

In Subsections B and C below, we show how the polytope expansion can be used in the computation of tensor product decompositions and subalgebra branching rules, respectively. In the first subsection, we show how weight multiplicities can be recovered easily from the polytope expansion coefficients (polytope multiplicities). The motivation is to demonstrate explicitly that the two descriptions of weight systems and characters are ultimately equivalent.

\subsection{Weight multiplicities}

To distinguish the {\it dominant} weight multiplicities, we define
\ben {\rm m}_\lambda(\mu)\ :=\ \left\{ \matrix{ {\rm mult}_\lambda(\mu)\,,\ & {\rm if}
\ \mu\in P_+\ ;\cr 0\ ,\ \ &\  {\rm otherwise}.
  }\right.\ \label{mlamu}\een
Let $W\mu$ denote the Weyl orbit of the weight $\mu$, and define an (even) Weyl orbit sum as \ben E_\lambda\ =\ \sum_{\varphi\in W\lambda}\, e^\varphi\ =\ \frac{\vert W\lambda\vert}{\vert W\vert}\sum_{w\in W}\, e^{w\lambda}\ .\label{Ela}\een
(These orbit sums have been  studied as functions on weight space in Ref.~\onlinecite{KPCS}.)
Then (\ref{Wsymm}) applied to (\ref{chmult}) gives \begin{equation}
{\rm ch}_\lambda\ =\ \sum_{\mu\leq \lambda}\, {\rm m}_\lambda(\mu)\, E_\mu\ ,
\label{chmE}\end{equation}
making the Weyl symmetry manifest. Clearly, determining the dominant weight multiplicities ${\rm m}_\lambda(\mu)$ is sufficient for knowledge of the ${\rm mult}_\lambda(\mu)$.

But the dominant weight multiplicities can be found easily from the polytope multiplicities, since
\ben B_\lambda\ =\ \sum_{\lambda\geq \mu\in P_+}\ E_\mu\ \ .\label{BE}\een   This implies \ben m_\lambda(\mu)\ =\ \sum_{\mu\leq\varphi\leq\lambda}\, {\rm polyt}_\lambda(\varphi)\ =:\ \sum_{\varphi=\mu}^{\lambda}\, {\rm polyt}_\lambda(\varphi)\ .\label{mA}\een

The expansions of Lie characters into Weyl-orbit sums  (\ref{Ela}) and into polytope sums (\ref{chAB}) are similar. Both the $E_\mu$ and the $B_\mu$ are multiplicity-1 sums. In general, however, the polytope sum $B_\mu$ contains many Weyl-orbit sums $E_\mu$, as (\ref{BE}) indicates.

The polytope multiplicities encode the weight system and multiplicities in a more economical manner than do the dominant weight multiplicities. Yet the latter can be recovered from the former simply (and non-negatively).

Let us emphasize here that the purpose of introducing the polytope expansion is not to provide another method of calculating weight multiplicities.  There are very efficient techniques known already for weight-multiplicity computations.  The message is rather that the polytope expansion and polytope multiplicities provide a different description of weight systems and characters, that carries the information in a more economical, but still useful way, and that also manifests patterns in weight systems beyond Weyl-group symmetry.

\subsection{Tensor product decompositions}
Taking the trace of the tensor product decomposition \ben L(\lambda)\otimes L(\mu)\ \hookrightarrow\ \bigoplus_{\nu\in P_+}\, T_{\lambda,\mu}^\nu\, L(\nu)\    \label{tpd}\een yields \ben  {\rm ch}_\lambda\, {\rm ch}_\mu\ =\ \sum_{\nu\in P_+}\, T_{\lambda,\mu}^\nu\,    {\rm ch}_\nu\ .  \label{chT}\een Substituting the Weyl character formula leads to a well-known formula  for the tensor product multiplicities: \ben T_{\lambda\mu}^\nu\ =\ \sum_{w\in W}\, (\det w)\,\, {\rm mult}_\mu(w.\nu-\lambda)\ \ .   \label{Tmult}\een This is equivalent to the Klimyk formula \cite{Kl} (or see Ref.~\onlinecite{Hum},  Exercise 24.9, for example), encoding the Racah-Speiser algorithm \cite{RS} for their computation. Knowledge of the multiplicities of weights, dominant and non-dominant, is necessary for this calculation.

The polytope expansion (\ref{chAB},\ref{polytA}) reduces that requirement to knowledge of the polytope multiplicities.  To take  advantage of the multiplicity-1 property of the polytope generating functions $B_\sigma$, we write \ben T_{\lambda,\mu}^\nu\ =\ \sum_{\sigma\in P_+}\, {\rm polyt}_\mu(\sigma)\, U_{\lambda,\sigma}^\nu \ \ ,\label{TAU}\een where
\ben {\rm ch}_\lambda\, B_\sigma\ =:\ \sum_{\nu\in P_+}\, U_{\lambda,\sigma}^\nu\,    {\rm ch}_\nu\ .\label{chBU}\een A Racah-Speiser calculation produces \ben U_{\lambda,\sigma}^\nu\ =\ \sum_{w\in W}\, (\det w)\,\,\delta_\sigma\,(w.\nu-\lambda)\ . \label{Udelta}\een Here \ben \delta_\sigma\,(\kappa)\ :=\ \left\{\matrix{ 1\, ,\ \ & \kappa\in Pt_\sigma\ ;\cr
0\, ,\ \ & \kappa\not\in Pt_\sigma\ .  }\right.  \label{deltaPt}\een

\subsection{Branching rules}
Suppose the semi-simple algebra $\bar X_{\bar r}$ is embedded in the algebra $X_r$, i.e. $\bar X_{\bar r}$ is a subalgebra of $X_r$. The branching rules \ben L(\lambda)\ \hookrightarrow\ \bigoplus_{\bar\mu\in \bar P_+}\, b_{\lambda,\bar\mu}\, \bar L(\bar\mu)\     \label{Lbr}\een result, and the branching coefficients $b_{\lambda, \bar\mu}$ are of interest. Here $\bar L(\bar\mu)$ indicates the $\bar X_{\bar r}$-representation of highest weight $\bar\mu$.  The latter is in $\bar P_+$, the set of dominant weights of the subalgebra $\bar X_{\bar r}$, i.e. the set of possible highest weights of irreducible $\bar X_{\bar r}$-representations. $\bar P \supset \bar P_+$ will indicate the set of $\bar X_{\bar r}$ weights, i.e. the weight lattice of $\bar X_{\bar r}$.

The ranks satisfy $\bar r\le r$, and the subalgebra implies that the weight space of $X_r$ projects onto the weight space of $\bar X_{\bar r}$. Let ${\bar I}$ denote the projection operator determined by the embedding, so that if $\nu\in P$, then ${\bar I} \nu\in \bar P$.  The branching coefficients $b_{\lambda, \bar\mu}$ of (\ref{Lbr}) also describe the decomposition of projected characters: \ben \bar I\,{\rm ch}_\lambda\ =\ \sum_{\bar\mu\in \bar P_+}\, b_{\lambda,\bar\mu}\, \bar{\rm ch}_(\bar\mu)\  .    \label{Ichbr}\een Here \ben   \bar I\,{\rm ch}_\lambda\ :=\ \sum_{\mu\in P}\, {\rm mult}_\lambda(\mu)\, e^{{\bar I} \mu}\ \ , \label{barch}\een if (\ref{chmult}) is obeyed.

Normally, the projection $\bar I$ is known, and the branching coefficients are desired.  They can be determined from (\ref{Ichbr}) by direct computation, but the procedure can be simplified in various ways. One way makes use of the character decompositions (\ref{chmE})  into orbit sums, and the inverse procedure.\cite{PS} We will describe it now, and then show that the polytope expansion of characters (or decomposition of characters into polytope functions) can be used in a completely analogous, but more efficient way.

For notational convenience, write the character-to-orbit decompositions for $X_r$ and $\bar X_{\bar r}$  as \ben {\rm ch}_\lambda\ =\ \sum_{\mu\leq \lambda}\, M_{\lambda,\mu}\, E_\mu\ , \qquad\&\qquad \bar{\rm ch}_{\bar\lambda}\ =\ \sum_{\bar\mu\leq \bar\lambda}\, {\bar M}_{\bar\lambda,\bar\mu}\, \bar E_{\bar\mu}\ ,    \label{MbarM}\een respectively. The matrices $(M_{\lambda,\mu})$ and ${\bar M}_{\bar\lambda,\bar\mu}$ have entries 1 on their diagonals and are triangular, so that they are easily inverted. That is, we can write \ben  E_\lambda\ =\ \sum_{\mu\leq \lambda}\, M^{-1}_{\lambda,\mu}\, {\rm ch}_\mu\ ,   \label{EMinvch}\een and similarly for the subalgebra $\bar X_{\bar r}$.   Furthermore,\cite{GJW} we have \ben    M^{-1}_{\lambda,\mu}\ =\  \frac{\vert  W\lambda\vert}{\vert W \vert}\,\, \sum_{x,w\in W}\, (\det w)\, \delta_{x\lambda+w.0,\, \mu+\rho}\   \qquad  ,\label{MinvW}\een and a similar formula for the subalgebra.

The procedure \cite{PS} is then to first decompose the algebra characters into orbit sums, then decompose the algebra orbit sums into subalgebra ones: \ben  \bar I\,E_\lambda\ =\ \sum_{\bar\mu}\, e_{\lambda, \bar\mu}\, {\bar E}_{\bar\mu}\ , \label{EobE}\een then reassemble the resulting orbit sums into subalgebra characters: \ben  b_{\lambda,\bar\lambda}\ =\ \sum_{\mu\leq\lambda}\, M_{\lambda,\mu}\, \sum_{\bar\mu\in P_+}\, e_{\mu,\bar\mu}\, \sum_{\bar\mu\geq \bar\lambda}\, {\bar M}^{-1}_{\bar\mu,\bar\lambda}\ .  \label{bMebMinv}\een Calculating the orbit-to-orbit branching coefficients is a relatively simple task, so this procedure is relatively efficient.

The point is simply that the same procedure can be applied when orbit sums are replaced by polytope sums $B_\lambda$, $\bar B_{\bar \lambda}$.  If the polytope-to-polytope branchings are\ben   \bar I\,B_\lambda\ =\ \sum_{\bar\mu}\, p_{\lambda, \bar\mu}\, {\bar B}_{\bar\mu}\  ,\label{BpbB}\een we replace $M_{\lambda,\mu}$ and ${\bar M}^{-1}_{\bar\mu,\bar\lambda}$ in (\ref{bMebMinv}) by the polytope multiplicities $A_{\lambda,\mu}$ and the inverses ${\bar A}^{-1}_{\bar\mu,\bar\lambda}$, to get \ben  b_{\lambda,\bar\lambda}\ =\ \sum_{\mu\leq\lambda}\, A_{\lambda,\mu}\, \sum_{\bar\mu\in P_+}\, p_{\mu,\bar\mu}\, \sum_{\bar\mu\geq \bar\lambda}\, {\bar A}^{-1}_{\bar\mu,\bar\lambda}\ .  \label{bApbAinv}\een The relation (\ref{AinvF}) plays the role of (\ref{MinvW}) here, so that the reassembly of polytope sums into characters can be done quite efficiently.




%
%

%

\begin{acknowledgments}
The author thanks Terry Gannon for reading the manuscript. This research was supported in part by a Discovery Grant from the Natural Sciences and Engineering Research Council (NSERC) of Canada.
\end{acknowledgments}



\end{document}